\documentclass[]{siamltex}
\usepackage{amsmath}
\usepackage{amssymb}
\usepackage{graphicx}

\usepackage{tikz}
\usepackage{pgfplots}
\usepackage[ComplainAboutUnexternalized]{externaltikz}
\usepackage{calc}

\usepackage{hyperref}

\newcommand{\R}{\mathbb{R}} 
\newcommand{\rhoK}{\rho_K} 

\newtheorem{remark}{Remark}

\begin{document}

\title{A nonlinear discrete-velocity relaxation model for traffic flow}

\author{R. Borsche\footnotemark[1] 
     \and  A. Klar\footnotemark[1] \footnotemark[2]}
\footnotetext[1]{Technische Universit\"at Kaiserslautern, Department of Mathematics, Erwin-Schr\"odinger-Stra{\ss}e, 67663 Kaiserslautern, Germany 
  (\{borsche, klar\}@mathematik.uni-kl.de)}
\footnotetext[2]{Fraunhofer ITWM, Fraunhoferplatz 1, 67663 Kaiserslautern, Germany} 
 
\date{}


\maketitle

\begin{abstract}
We derive  a nonlinear 2-equation discrete-velocity model for traffic flow from a continuous kinetic model. The model converges to  scalar Lighthill-Whitham type equations in the relaxation limit for all ranges of traffic data. Moreover, the model has an invariant domain  appropriate for traffic flow modeling.
It  shows some similarities with the Aw-Rascle traffic model. However, the new model is simpler and yields, in case of a concave fundamental diagram, an example for a  totally linear degenerate hyperbolic relaxation model.
We discuss the details of the hyperbolic main part and  consider boundary conditions for the limit equations derived from the relaxation model. Moreover, we  
investigate the cluster dynamics of the model for vanishing braking distance  and consider a relaxation scheme build on the kinetic discrete velocity model. Finally, numerical results for various situations are presented, illustrating the analytical results.
\end{abstract}




\section{Introduction}

Starting with the work of Lighthill and Whitham \cite{Whi74}, there have been many approaches to a continuous modeling of traffic flow problems.
Macroscopic models are usually based on scalar hyperbolic equations like the above cited model or systems of hyperbolic equations with relaxation term, see \cite{Pay79} for a classical equation. 
More recently, an improved traffic flow model using hyperbolic systems with relaxation has been presented by Aw and Rascle \cite{AR}.
For discussions and extensions see, for example, \cite{AKMR,Deg,Ber,G,R}.
On the other hand, kinetic equations have also been widely used as a tool to model traffic flow problems, see \cite{PH71,PF75, Nel95,Hel95B, KW981}.
In \cite{KW97,KW00} non-local terms are introduced into the equations to guarantee information transport against the flow direction. 
We refer to \cite{BD11} for a recent mathematically oriented review and further references.
In a simplified context after discretizing the velocity space and neglecting non-localities in the kinetic model, the resulting discrete velocity models are hyperbolic systems 
with relaxation terms.
These so called relaxation systems have been widely used for example for numerical purposes, see \cite{JX95}.

However, a naive application of relaxation systems  in the case of traffic flow leads to similar problems as for the full kinetic equation.
Either negative discrete velocities are allowed, which is not meaningful from the traffic flow point of view, or
some kind of non-locality has to be introduced into the equations, see the next section for a detailed discussion. 
For a non-local discrete velocity traffic model, we refer to \cite{HPS2}.
From the point of view of hyperbolic relaxation systems this  is closely related to  the so-called subcharacteristic condition \cite{CL93}.
A further complication is given by the fact, that for traffic flow modeling the state space is restricted to positive and bounded velocities and densities. 
This leads to requirements on the invariant domains of the equations, see  again the next section for details.
Considering discrete velocity relaxation models in the usual form without explicit non-localities there is no way to achieve a correct invariant domain with a linear hyperbolic part of the relaxation system.

In this paper we aim at deriving and investigating  a discrete velocity model with nonlinear hyperbolic part fulfilling the above requirements.
In particular, we require that the model has the correct invariant domains and fulfills the sub-characteristic condition and converges to  
a  scalar Lighthill-Whitham type equation in the relaxation limit.
It will turn out that the resulting model has some similarities with the AW-Rascle model,
being a  hyperbolic model of the so called Temple class and in a special, but relevant, case being an example for a  totally linear degenerate hyperbolic equation with relaxation term.

The paper is organized in the following way. 
In section \ref{notations} we discuss classical discrete velocity relaxation models and their drawbacks in the traffic flow case. Moreover, we discuss the relation to the Aw-Rascle model and its modifications.
In section \ref{discretevelocity}  we derive a new nonlinear discrete velocity kinetic model for traffic flow from a continuous kinetic traffic equation.
We consider the associated macroscopic equations and the convergence to the Lighthill-Whitham equations.
In the subsequent section \ref{macroscopicmodel} the macroscopic equations are discussed in detail including hyperbolicity, integral- and shock curves and Riemann invariants of the homogeneous system.
In section \ref{boundaryconditions} we consider the derivation of boundary conditions for the limiting Lighthill-Whitham type equations from the boundary conditions of the underlying kinetic problem based on the analysis of the kinetic boundary layer. In section \ref{relaxationmethods} we discuss a relaxation method based on the nonlinear discrete velocity model.
Section \ref{cluster} discusses a constrained linear model derived from the kinetic  model in the  limit of small braking distances. 
Finally,
numerical results are presented in Section \ref{numericalresults}.

\section{Notations and motivation}
\label{notations}
The most important tool in traffic flow modeling is the fundamental diagram $F(\rho) , 0 \le \rho \le 1$.
We consider smooth  functions $F$ with $F(0) = F(1) =0$ and the following property. There is a $0 < \rho^{\star} < 1$ such that 
$F^{\prime}(\rho) >0$ for $ 0\le \rho < \rho^{\star} $ and $F^{\prime}(\rho) <0$ for $ \rho^{\star}< \rho < 1 $.
We use the notation $\tau (\rho)$ for the value $\tau(\rho) \neq \rho$ such that $F(\tau(\rho)) = F(\rho)$, compare \cite{CP02}.
Here and in the remainder of the paper we set the maximal density to $1$ as well as the maximal velocity.

Discrete velocity models have been investigated in many works, see
\cite{Illner} for a review. 
To begin with, we consider a classical discrete velocity model for the distribution functions $f_1$ and $f_2$ associated to the two velocities $-1 \le v_1 \le v_2 \le 1$. The equilibrium functions are $0 \le f^e_1 (\rho),f^e_2(\rho) \le 1$
where the density $\rho$ is given by $\rho = f_1 + f_2 =f^e_1 + f^e_2$ and the mean velocity by $q = v_1 f_1 + v_2f_2$. Moreover, the equilibrium flux is  $v_1 f^e_1 + v_2 f^e_2 = F(\rho)$. The equations are
\begin{eqnarray}
\label{dvm0}
\partial_t f_1 + v_1 \partial_x f_1 = -\frac{1}{\epsilon} \left(f_1-f^e_1(\rho) \right)  \\
\partial_t f_2 + v_2 \partial_x f_2 =  -\frac{1}{\epsilon} \left(f_2-f^e_2 (\rho)\right).
\end{eqnarray}
From the above we have
\begin{align*}
f^e_1 = \frac{v_2 \rho-F(\rho)}{v_2-v_1} \  ,\ 
f^e_2= \frac{ F(\rho)-v_1 \rho}{v_2-v_1}
\end{align*}
and 
$$
f_1 =  \frac{v_2 \rho-q}{v_2 -v_1}\ ,\ 
f_2= \frac{q-v_1 \rho }{v_2 -v_1}\ .
$$
Equation (\ref{dvm0}) can be rewritten as
\begin{align}\label{dvm}
\begin{aligned}
\partial_t \rho +  \partial_x q &=0\\
\partial_t q   -v_1 v_2 \partial_x \rho + (v_1+v_2)\partial_x q  & =-\frac{1}{\epsilon} \left(q-F(\rho) \right) \ .
\end{aligned}
\end{align}

We remark first  that the invariant region of the equations is given by the rectangle $0 \le f_1, f_2 \le 1$, which gives $0\le v_2 \rho -q \le v_2-v_1$ and $0\le q-v_1 \rho  \le v_2-v_1$. The latter is rewritten as
$v_2\rho - v_2-v_1\le q \le v_2 \rho  $ and $v_1\rho  \le q \le v_1 \rho+v_2-v_1$.

However, for a reasonable discrete velocity traffic model the  invariant region 
should  be given by the triangle $ 0 \le \rho \le 1$ and  $0 \le q \le \rho$ or in terms of $f_1,f_2$ by  the region $0\le f_1 \le 1$ and 
$0 \le f_1+f_2 \le 1$. 
We note that $q \le \rho$ is a bound for 
the maximal velocity. See the discussion in \cite{AR,AKMR,Deg}.

Having discussed this, we remark that 
the region $ 0 \le \rho \le 1, 0 \le q \le \rho$ is not an invariant domain of the above equations.
For the above discrete velocity model there is no guarantee for positive  $q$ except for the case $v_1,v_2 \ge 0$, but in this case the bound $\rho \le 1$ is not satisfied. 
Indeed, we observe that for   $v_1=0$ and $v_2=1$  the triangle $ 0 \le \rho \le 1$, $0 \le q \le \rho$  is contained in the invariant domain $\rho - 1\le q \le  \rho  $ and $0  \le q \le 1$,
but is still not an invariant domain itself.

Moreover, under the  subcharacteristic condition  $v_1 \le F^\prime (\rho) \le v_2$ we obtain convergence of (\ref{dvm0}) to the conservation law
$$
\partial_t \rho + \partial_x F(\rho) =0\ .
$$

Obviously, $v_1 \ge 0$ does not allow for  situations where $F^\prime (\rho )$ is negative. In particular, situations with traffic jams are not treated.
We remark, that one cannot remedy this by manipulating the right hand side as in \cite{HPS}, where an unstable relaxation system has been used.

In the present paper we introduce a new 
discrete-velocity model for traffic flow being on the one hand a reasonable model for traffic flow in the sense that  with $v_1,v_2 \ge 0$ one obtains
an invariant region given by 
$0 \le f_1, f_2 \le 1$ and $0 \le f_1+f_2 \le 1$, or equivalently $ 0 \le \rho \le 1, 0 \le q \le \rho$, and 
being on the other hand a model which  converges as $\epsilon$ goes to $0$ to the conservation law for all values of $\rho$.

We remark that the modified Aw-Rascle equations \cite{Deg} fulfill the above requirements. The relaxation property has been investigated for example in  \cite{AKMR,R}. In fact, as will be seen later, our new model shares  some of the properties of  the modified Aw-Rascle model.
However, the Aw-Rascle model  is not derived from a kinetic discrete-velocity model. Moreover,  the new model is different from a modeling point of view. In the Rascle model acceleration and breaking influence the hyperbolic part of the equations as well as the relaxation term. 
 Our model uses  only the non-locality in the breaking term  for a  contribution to the hyperbolic part of the equations. All other physical influences are summarized in the relaxation term.
  We refer to  the original work on kinetic traffic flow equations \cite{PH71} for similar considerations. In this sense the model is a minimal model using only those physical phenomena in the hyperbolic part of the equations which are necessary to guarantee the subcharacteristic condition and convergence to the scalar conservation law.

\section{Discrete velocity traffic model}
\label{discretevelocity}
Our starting point is a kinetic equation with continuous velocity space, compare  \cite{KW97,KW981,KW00} and  see also \cite{Hel95B,PF75,PH71,Nel95}. For $t \in \R^+$, $x \in \R$, $v \in[0,1]$ and the distribution function $f=f(x,v,t)$ we consider
the equation
$$
\partial_t f + v \partial_x f = J_R (f) + J_{NL} (f)
$$
with a relaxation term $J_R$ relaxing to an equilibrium function $f_0(\rho)$
\begin{align*}
J_R (f)=  -\frac{1}{\epsilon} \left( f - f_0 (\rho) \right)\ ,
\end{align*}
with $\int f_0 (v) dv = \rho $  and $\int v f_0(v) dv = F(\rho)$.
Additionally we consider a term containing the  non-local effects due to braking interactions $J_{NL}(f)$ given by
$$
J_{NL} (f) = J_B (f,H) - J_B (f,0)\ ,
$$
where the braking term $J_B (f,H)$ is given as in  \cite{BKK}. We obtain in a simplified case
\begin{align*}
 J_B (f,H) = \frac{1}{1-\rho }\int_{\hat v > v}  (\hat v -v)  f(x,\hat v) f(x+H, v ) d \hat v\\
- \frac{1}{1-\rho }\int_{\hat v < v } ( v - \hat v)  f(x, v) f(x+H, \hat v ) d \hat v\ .
\end{align*}
Here, $H$ is a measure for the  look-ahead and  the non-locality of the equations. 
The underlying microscopic model contains a braking interaction, where the driver at $x$ with velocity $v$ reacts to his predecessor at $x+H$ with velocity $\hat v $  if 
$\hat v < v$. The new velocity resulting of this braking interaction is exactly the velocity of the leading car.
Moreover, the interaction strength is modulated by a factor $\frac{1}{1-\rho}$ to increase the frequency of breaking interaction for dense traffic.
This can be approximated by
\begin{align*}
 J_B (f,H) 
\sim \frac{1}{1-\rho } \int_{\hat v > v}  ( \hat v -v )  f(\hat v) (f( v ) + H \partial_x f (v))d \hat v \\- \frac{1}{1-\rho }\int_{\hat v < v}  ( v-\hat v  ) f( v) (f( \hat v ) + H \partial_x f (\hat v))  d \hat v\ 
\end{align*}
 and thus
\begin{align*} 
J_{NL} 
\sim \frac{H}{1-\rho } \int_{\hat v > v}  ( \hat v -v )  f(\hat v)   \partial_x f (v)d \hat v- \frac{H}{1-\rho } \int_{\hat v < v}  ( v-\hat v  ) f( v)   \partial_x f (\hat v)  d \hat v\ .
\end{align*}

For  a two velocity model with the velocities $0 \le v_1<v_2 \le 1$ 
we obtain, setting 
 $v=v_1$ and $v=v_2$ respectively,  the same relaxation  term as in equation (\ref{dvm0}). The non-local term yields for $v=v_1$
\begin{eqnarray*}
J_{NL,1} =  \frac{H}{1-\rho }  ( v_2 -v_1)  f_2\partial_x f_1
\end{eqnarray*}
and for 
$v = v_2$
\begin{eqnarray*}
J_{NL,2} = -\frac{H}{1-\rho }  ( v_2 -v_1) f_2 \partial_x f_1\ .
\end{eqnarray*}
Alltogether, the two velocity discrete-velocity model is given by the following nonlinear discrete velocity model
\begin{eqnarray*}
\partial_t f_1 + v_1 \partial_x f_1 - \frac{H}{1-\rho }  ( v_2 -v_1)  f_2\partial_x f_1= -\frac{1}{\epsilon} \left(f_1-\frac{v_2 \rho-F(\rho)}{v_2-v_1} \right)  \\
\partial_t f_2 + v_2 \partial_x f_2 + \frac{H}{1-\rho }  ( v_2 -v_1)  f_2\partial_x f_1=  -\frac{1}{\epsilon} \left(f_2-\frac{ F(\rho)-v_1 \rho}{v_2-v_1}\right).
\end{eqnarray*}
Using
$
f_1 = \frac{v_2 \rho-q}{v_2 -v_1}
$
and
$
f_2= \frac{q-v_1 \rho }{v_2 -v_1}
$
gives the macroscopic equations for density $\rho$ and mean flux $q$
\begin{align}\label{macro}
\begin{aligned}
\partial_t \rho + \partial_x q &=0\\
\partial_t q + \partial_x  P(\rho,q ) - (v_1 + v_2)\frac{H}{1-\rho}   \frac{q-v_1 \rho}{v_2-v_1} \left(v_2 \partial_x \rho-\partial_x q\right) & =-\frac{1}{\epsilon} \left(q-F(\rho) \right) \ 
\end{aligned}
\end{align}
with
$
P(\rho,q) = (v_1+v_2) q-v_1 v_2 \rho.
$
In the following we consider the canonical choice of  the two velocities as $v_1=0,v_2 =1$. Otherwise, situations with very low or very high velocites could not be covered.
With  $v_1=0,v_2 =1$  we get
$$
f_1 =\rho-q, f_2 =q
$$
and the discrete velocity model
\begin{eqnarray*}
\partial_t f_1  - \frac{H}{1-\rho }    f_2\partial_x f_1= -\frac{1}{\epsilon} \left(f_1-\rho +F(\rho) \right)  \\
\partial_t f_2  + \partial_x f_2 + \frac{H}{1-\rho }   f_2\partial_x f_1=  -\frac{1}{\epsilon} \left(f_2- F(\rho) \right)
\end{eqnarray*}
or the macroscopic equation
\begin{align}\label{macro0}
\begin{aligned}
\partial_t \rho + \partial_x q &=0\\
\partial_t q + \frac{H q}{1-\rho}  \partial_x \rho  + (1-\frac{H q}{1-\rho})\partial_x  q  &=-\frac{1}{\epsilon} \left(q-F(\rho) \right) \ .
\end{aligned}
\end{align}

The details of the hyperbolic equation will be considered in the next section. 
Concerning the convergence of the equations towards the scalar conservation law $\partial_t \rho + \partial_x F(\rho) =0$ as
$\epsilon $ tends to $0$ we have to assure the subcharacteristic condition.
The eigenvalues of the system are 
$$
\lambda_1 = - \frac{H q}{1-\rho}\ , \ \lambda_2 =1\ .
$$
Setting $q = F(\rho)$ we get
$$
\lambda_1  = - \frac{ H F(\rho)}{1-\rho}\ ,\ \lambda_2 =1\ .
$$
The subcharacteristic condition gives 
$$
-  \frac{ H F(\rho)}{1-\rho} \le F^\prime(\rho) \le 1 \  \mbox{ for } \  0 \le \rho \le 1\ .
$$

\begin{remark}
In the classical Lighthill Whitham case with $F(\rho) =\rho (1-\rho) $ and  $F^\prime (\rho) = 1- 2 \rho$ this yields the condition
$$
- H  \rho \le 1- 2 \rho \le 1\  \mbox{ for } \ 0 \le \rho \le 1\ .
$$
This is satisfied for all $H\ge 1$.
\end{remark}

\begin{remark}
	In the case of a general traffic fundamental diagram  $F(\rho)  $ the above condition is on the one hand guaranteeing that $F(\rho) \le \rho$, which guarantees that the equilibrium function lies in the invariant domain, see the next section. On the other hand, the first inequality yields 
	$$
	 F(\rho) \ge \frac{F(\rho^{\star})}{ (1-\rho^{\star})^H} (1-\rho)^H.
	$$
	 That means the value of $H$ must be  chosen according to the behaviour of the fundamental diagram at $\rho=1$.
	 In case of a concave flux function $F$  it  is sufficient to choose $H=1$.
\end{remark}

\begin{remark}
The strength of the breaking interaction can be also changed by changing the  coefficient 
 $\frac{1}{(1-\rho)}$ in front of the braking term into $ \frac{1}{(1-\rho)^n}, n >1$.
\end{remark}

\section{The nonlinear macroscopic  system}
\label{macroscopicmodel}
In this section we consider the macroscopic 
 hyperbolic system  in more detail.
First, we consider the homogeneous system
\begin{align}\label{macrohom}
\begin{aligned}
\partial_t \rho + \partial_x q &=0\\
\partial_t q + \frac{H q}{1-\rho}  \partial_x \rho  + (1-\frac{H q}{1-\rho})\partial_x  q  &=0 \ .
\end{aligned}
\end{align}
The eigenvalues of the system are 
$$
\lambda_1  = - \alpha = - \frac{H q}{1-\rho} < \lambda_2 =1
$$
with eigenvectors
\begin{equation}\label{eigenvec}
r_1 = \left(
\begin{array}{c}
1\\
- \alpha
\end{array}
\right) ,
r_2 = \left(
\begin{array}{c}
1\\1
\end{array}
\right)\ .
\end{equation}
Moreover, for the characteristic families we obtain  
\begin{align*}
\nabla \lambda_1 \cdot r_1 = 
\left(
\begin{array}{c}
-\frac{Hq}{(1-\rho)^2}\\
- \frac{H}{1-\rho}
\end{array}
\right)  \cdot
\left(
\begin{array}{c}
1\\- \frac{Hq}{1-\rho}
\end{array}
\right)\ 
= \frac{Hq}{(1-\rho)^2}(H-1)
\end{align*}
and 
\begin{align}
\nabla \lambda_2 \cdot r_2 = 0\ .
\end{align}
That means the $r_1$-field is genuinely nonlinear for $H\neq 1$ and $q \neq 0$, the $r_2$-field is linearly degenerate.
For $H=1$ we have a (totally) linear degenerate system.

The integral curves of the system are determined 
considering  the following ODE's.
For the 1-field we have
\begin{align*}
\begin{aligned}
\rho^\prime = 1\ ,&&
q^\prime =- \frac{Hq}{1-\rho}\ .
\end{aligned}
\end{align*}
This gives the integral curve $q= q_L \frac{(1-\rho)^H}{(1-\rho_L)^H}$.

For the 2 field one obtains
\begin{align*}
\begin{aligned}
\rho^\prime = 1\ ,&&
q^\prime =1\ .
\end{aligned}
\end{align*}
That means the 2-integral curves are straight lines $q = \rho -\rho_R+q_R$  with slope $1$.

For $H=1$ the dynamics are completely described by the integral curves.
In case $H>1$ additionally  the shock curves have to be investigated for the 1-field.

We rewrite the equations in
conservative form.
For general $H$, we choose the variable $$z =  \frac{Hq}{(1-\rho)^H}\ .$$ 
Note that for $H=1$ we have $z= - \lambda_1$ as expected.
Then
\begin{align*}
\begin{aligned}
\partial_t z 
=& \frac{H^2q}{(1-\rho)^{H+1}} \partial_t \rho + \frac{H}{(1-\rho)^H} \partial_t q 
= \frac{H}{(1-\rho)^H} \partial_t q -\frac{H^2q}{(1-\rho)^{H+1}} \partial_x q \\
=&  -\frac{H}{(1-\rho)^H} \left(\frac{H q}{1-\rho}  \partial_x \rho  + (1-\frac{H q}{1-\rho})\partial_x  q +\frac{1}{\epsilon} \left(q-F(\rho) \right)\right)-\frac{H^2q}{(1-\rho)^{H+1}} \partial_x q \\
=&  - \frac{H^2 q}{(1-\rho)^{H+1}} \partial_x \rho  + \left( -\frac{H }{(1-\rho)^H} + \frac{H^2 q }{(1-\rho)^{H+1}} - \frac{H^2 q}{(1-\rho)^{H+1}} \right)\partial_x  q\\
&-\frac{1}{\epsilon}\frac{H \left(q-F(\rho) \right)}{(1-\rho)^H} 
\ .\end{aligned}
\end{align*}
We obtain
\begin{align*}
\begin{aligned}
\partial_t z &=  - \frac{H^2 q}{(1-\rho)^{H+1}} \partial_x \rho  -  \frac{H }{(1-\rho)^H} \partial_x  q -\frac{1}{\epsilon}\frac{H \left(q-F(\rho) \right)}{(1-\rho)^H}\\
&=  -\partial_x z-\frac{1}{\epsilon}\frac{H \left(q-F(\rho) \right)}{(1-\rho)^H}
\end{aligned}
\end{align*}
and the conservative system
\begin{align}\label{cons}
\begin{aligned}
\partial_t \rho + \partial_x q&=0\\
\partial_t z + \partial_x z &= -\frac{1}{\epsilon}\frac{H }{(1-\rho)^H}\left(q-F(\rho) \right)
\ .\end{aligned}
\end{align}
In closed form this is 
\begin{align}\label{consclosed}
\begin{aligned}
\partial_t \rho + \frac{1}{H}\partial_x(z (1-\rho)^H)&=0\\
\partial_t z + \partial_x z &= -\frac{1}{\epsilon}\left( z -\frac{H F(\rho) }{(1-\rho)^H}  \right)\ . \end{aligned}
\end{align}

The eigenvalues are written as 
$\lambda_1 = -z(1-\rho)^{H-1}$ and $\lambda_2 = 1$.
The eigenvectors in conservative variables are
\begin{equation*}
r_1 = \left(
\begin{array}{c}
1\\
0
\end{array}
\right)\ ,\ 
r_2 = \left(
\begin{array}{c}
1\\ \frac{zH}{1-\rho}  + \frac{H}{(1-\rho)^H}
\end{array}
\right)\ .
\end{equation*}

The integral curves are in conservative form for the 
1-field given by straight lines  $z=z_L = const$.
The integral curves for the 
2-field are given by

\begin{align}\label{incurve}
\begin{aligned}
z(\rho)  =  H\frac{\rho-\rho_R}{(1-\rho)^H}+ \frac{z_R(1-\rho_R)^H}{(1-\rho)^H}\ .
\end{aligned}
\end{align}
For the 1-field we have to consider additionally the 
1-shock curves.
The Rankine-Hugoniot conditions give
\begin{align}\label{locus}
\begin{aligned}
\frac{1}{H}z(1-\rho)^H - \frac{1}{H}z_{L/R}(1-\rho_{L/R})^H = s (\rho-\rho_{L/R})\\
z- z_{L/R} = s(z-z_{L/R})\ .
\end{aligned}
\end{align}

This yields either $z=z_L$. 
That means also for the  1-field shocks and integral curves coincide.
The speed of the shock (living on the shock curve with $\rho > \rho_L$) is computed from 
\begin{align}\label{s}
\begin{aligned}
z_L(1-\rho_L)^H - z_{R}(1-\rho_{R})^H &= s H (\rho_L-\rho_{R})\\
z_L&= z_{R}\  .
\end{aligned}
\end{align}
This gives
\begin{align}\label{shockspeed}
\begin{aligned}
s= \frac{z_L}{H}\frac{(1-\rho_L)^H - (1-\rho_{R})^H}{\rho_L-\rho_{R}}\ .
\end{aligned}
\end{align}

The second solution of (\ref{locus}) is $s=1$. This is the velocity of the 2-waves.
The 2-shock curve must be  the same as the two integral curve since the 2-field is linear degenerate.
Indeed,the 2-shock curve is given by 

\begin{align*}
\begin{aligned}
z(1-\rho)^H - z_R(1-\rho_R)^H =  H (\rho-\rho_R)
\end{aligned}
\end{align*}
or
\begin{align*}
\begin{aligned}
z(\rho)  =  H\frac{\rho-\rho_R}{(1-\rho)^H}+ \frac{z_R(1-\rho_R)^H}{(1-\rho)^H}
\end{aligned}
\end{align*}
as before 
 the 2-integral curves. Figure \ref{state} shows the integral curves in $(\rho,q)$ and $(\rho,z)$ variables for $H=1$.

 \begin{figure}[h]
 	\externaltikz{Statespace_H1}{
 	\begin{tikzpicture}[scale = 4]
 	\def\dr{0.2}
 	\node[below] at (0,0) {$(0,0)$};
 	\node[below] at (1,0) {$(1,0)$};
 	\node[right] at (1,1) {$(1,1)$};
 	\draw[->](0,0)--(1.2,0) node[below]{$\rho$};
 	\draw[->](0,0)--(0,1.2) node[left]{$q$};
 	\draw[->](0.25,0.54)--(0.63,0.54) node[above]  at (0.25,0.54){$1$-curve};
 	 	\draw[->](0.5,0.8)--(0.78,0.58) node[above] at (0.4,0.8){$2$-curve};
 	\draw(1,0)--(1,1);
 	\foreach \r in {0, \dr,...,1}
 	\draw(\r,0)--(1,1-\r);
 	\foreach \r in {0, \dr,...,1}
 	\draw(\r,\r)--(1,0);
 	\end{tikzpicture}
 		\begin{tikzpicture}[scale = 4]
 	 	\def\dr{0.1}
 	 	\def\dz{0.2}
 	 	\node[below] at (0,0) {$(0,0)$};
 	 	\node[below] at (1,0) {$(1,0)$};
 	 	\draw[->](0,0)--(1.2,0) node[below]{$\rho$};
 	 	\draw[->](0,0)--(0,1.2) node[left]{$z$};
 	 	\draw(1,0)--(1,1+\dz);
 	 	\foreach \z in { 0,\dz,...,1}
 	 	\draw[domain=\z/(1+\z):1,smooth,variable=\x,] plot ({\x},{\z});
 	 	\foreach \r in {0.5,0.6,...,1}
 	 	\draw[domain=0.0:1+\dz,smooth,variable=\x,] plot ({max(0,(\x+2*\r-1)/(\x+1))},{\x});
 	 	\draw[->](0.23,0.55)--(0.61,0.55) node[above]  at (0.23,0.55){$2$-curve};
 	 	 	\draw[->](0.3,0.8)--(0.56,0.60) node[above] at (0.3,0.8){$1$-curve};
 	 	\end{tikzpicture}
 	 }
 	\caption{State space and Riemann invariants in $(\rho,q)$ and $(\rho,z)$ variables for $H=1$.}
 	\label{state}
 \end{figure}
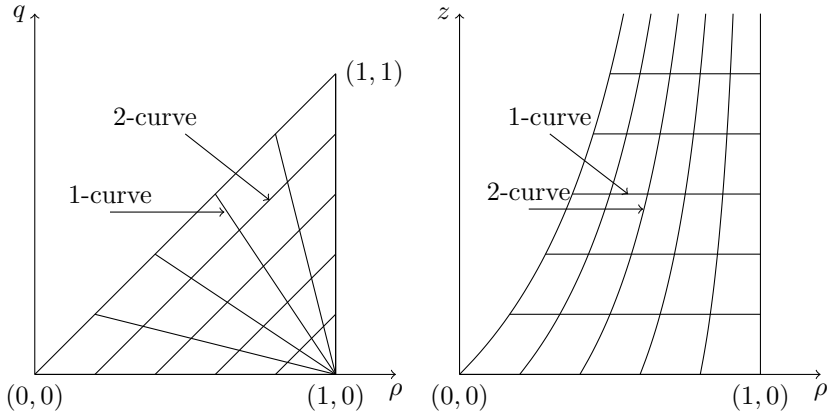
 
 \begin{figure}[h] 
 \centering
 \externaltikz{StatespaceH2}{
 	\begin{tikzpicture}[scale = 4]
 	\def\dr{0.2}
 	\node[below] at (0,0) {$(0,0)$};
 	\node[below] at (1,0) {$(1,0)$};
 	\node[right] at (1,1) {$(1,1)$};
 	\draw[->](0,0)--(1.2,0) node[below]{$\rho$};
 	\draw[->](0,0)--(0,1.2) node[left]{$q$};
 	\draw(1,0)--(1,1);
 	\foreach \r in {0, \dr,...,1}
 	\draw(\r,0)--(1,1-\r);
 	\foreach \r in {0, \dr,...,0.9}
 \draw[domain=\r:1,smooth,variable=\x,] plot ({\x},{\r*(1-\x)^2/(1-\r)^2});
 	\end{tikzpicture}
 	\begin{tikzpicture}[scale = 4]
 	 	 	\def\dr{0.1}
 	 	 	\def\dz{0.2}
 	 	 	\node[below] at (0,0) {$(0,0)$};
 	 	 	\node[below] at (1,0) {$(1,0)$};
 	 	 	\draw[->](0,0)--(1.2,0) node[below]{$\rho$};
 	 	 	\draw[->](0,0)--(0,1.2) node[left]{$z$};
 	 	 	\draw(1,0)--(1,1+\dz);
 	 	 	\foreach \z in {0.2,0.4,...,1}
 	 	 	\draw[domain={(1+\z-sqrt(1+2*\z))/\z}:1,smooth,variable=\x,] plot ({\x},{\z});
 	 	 	\foreach \r in {0.27,0.4,0.6,0.8,1}
 	 	 	\draw[domain=0.01:1+\dz,smooth,variable=\x,] plot ({max(0,(1+\x-sqrt(1+3*\x-4*\x*\r+\x*\r^2))/(\x))},{\x});
 	 	 	\end{tikzpicture}	
 	 	 } 	
 	\caption{State space and Riemann invariants in $(\rho,q)$   and $(\rho,z)$  variables for $H=2$.}
 	\label{state}
 \end{figure}
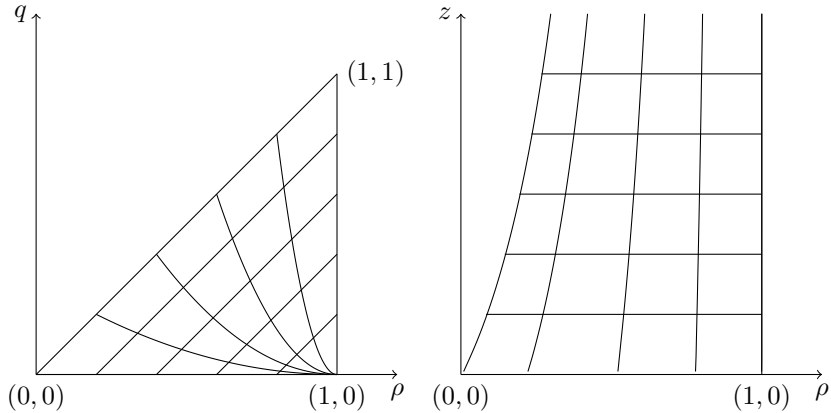
 
Note that the 2-field is linear degenerate and the 1-field has 
in conservative variables straight lines as integrals curves. Such systems have been investigated in several works, see, for example,
 \cite{Tem}.
Compare also  the structure of the  original Aw-Rascle traffic system \cite{AR}, which has the same property.

Moreover, as mentioned above, in the special case $H=1$ the system is even totally linear degenerate.

The  two Riemann invariants are easily determined:
the 1-Riemann invariant is $f_1 = \rho-q$, the 2-Riemann invariant is $\frac{H q}{(1-\rho)^H} $.

Finally, we observe that the integral curves for the 2-field are the same for any $H$. Thus,   it is easy to see that the region 
$0 \le \rho \le 1$, $0 \le q \le \rho$ is an invariant region for  the system  for all $H >0$.
This means for $u = \frac{q}{\rho}$ we have the invariant region $0 \le \rho \le 1, 0 \le u \le 1$
as it should be for a reasonable traffic flow equation.

\begin{remark}
One might compare the above equations to the  modified Aw-Rascle equations \cite{Deg}.
The second equation of the Aw-Rascle model in conservative form for the variable  $y = q +  \rho P(\rho)$ is given by
\begin{align}\label{macro0b}
\begin{aligned}
\partial_t y +   \partial_x \frac{q}{\rho} y    =-\frac{1}{\epsilon} \left( q-  F(\rho)\right) \ .
\end{aligned}
\end{align}
This can be rewritten as
\begin{align}\label{macrorascle}
\begin{aligned}
\partial_t q +   \tilde \alpha \partial_x \rho + \frac{q}{\rho} (1- \tilde \alpha) \partial_x q   =-\frac{1}{\epsilon} \left( q-  F(\rho)\right) \ .
\end{aligned}
\end{align}
with $ \tilde \alpha = q P^\prime (\rho)- \frac{q^2}{\rho^2}$, which shows some similarities with the equations considered here.
We remark that the second eigenvalue of the Aw-Rascle system 
 $\frac{q}{\rho}- \rho P^\prime (\rho)$, the first one is $\frac{q}{\rho}$ ,  does not have a fixed sign
compared to the present model. 
\end{remark}

\begin{remark}{Linearized relaxation system.}
The linearization of the above equations is
\begin{align}\label{lin}
\begin{aligned}
\partial_t \rho + \partial_x q &=0\\
\partial_t q + \alpha   \partial_x \rho  + (1-\alpha)\partial_x  q & =-\frac{1}{\epsilon} \left(q-F(\rho) \right) \ ,
\end{aligned}
\end{align}
where $\alpha \ge 0$.

The
characteristic variables of the linearized system are 
\begin{align*}
g_1 &=    \rho-q =f_1\\
g_2 &=   \alpha \rho+q  =  \alpha f_1+  (1+\alpha) f_2\ .
\end{align*}
This yields
\begin{eqnarray*}
\partial_t g_1 - \alpha \partial_x g_1 &= -\frac{1}{\epsilon} \left(-\frac{g_2-\alpha g_1}{1+\alpha}+F(\rho) \right)  \\
\partial_t g_2 + \partial_x g_2&=  -\frac{1}{\epsilon} \left(\frac{g_2-\alpha g_1}{1+\alpha}- F(\rho) \right)
\ .
\end{eqnarray*}
Thus, from the point of view of boundary conditions we have to  specify $\alpha f_1 + (1+\alpha) f_2$ at the left boundary and  $f_1$
 on the right boundary.
We note that the linearized equations face the same problems concerning their invariant domains as the naive relaxation model mentioned in Section \ref{notations}.

\end{remark}

In the  next step we investigate the boundary value problem for the nonlinear  kinetic system and consider the resulting 
boundary conditions for the limiting scalar hyperbolic problem as $\epsilon $ goes to $0$.

\section{Boundary conditions for the macroscopic equations derived from nonlinear kinetic equation}
\label{boundaryconditions}
In this section we determine boundary conditions for the scalar conservation law  from the  boundary value problem of the nonlinear kinetic relaxation system.
The boundary conditions for the limit equation  are obtained from the kinetic boundary conditions
considering a kinetic half-space problem at the boundary. We refer to  \cite{BLP79,CGS,BSS84} for boundary layers of kinetic equations  and to \cite{NT01,N96,WY99,LX96,LY01,WX99,AM04} for investigations of boundary layers for hyperbolic relaxation systems. 
 
The general procedure is as follows: the half space problem is determined by   a rescaling $x \rightarrow \frac{x}{\epsilon}$ of the spatial coordinate in the boundary layer. The boundary condition for the layer problem is given by the original kinetic boundary condition. The 
boundary condition for the limit equation is found from the asymptotic value of the half-space problem at infinity.

In the following we investigate first the kinetic layer equations and their asymptotic states and then  use the results to determine the boundary conditions for the macroscopic problem.

\subsection{Layer solution for nonlinear equations}
Let the left boundary be located at $x_L$.
Starting from equation (\ref{macro0}) and rescaling space as $x \rightarrow \frac{x-x_L}{\epsilon}$ one obtains 
the layer equations  for the left boundary  for $x \in [0, \infty)$ as 
\begin{align*}
\begin{aligned}
\partial_x q &=0\\
 \frac{Hq}{1-\rho}  \partial_x \rho  + (1-\frac{Hq}{1-\rho})\partial_x  q & =- \left(q-F(\rho) \right) \ 
\end{aligned}
\end{align*}
or
\begin{align}\label{layer}
\begin{aligned}
 q &=C\\
 \partial_x \rho   &= (1-\rho) \frac{F(\rho)-C}{HC}\ .
\end{aligned}
\end{align}
For $C < F(\rho^{\star})$
the above problem has three  fix points
$$\rho_{1} \le \rho^{\star}\ ,\ 
\rho_2 = \tau (\rho_1) \ge \rho^{\star}\ ,\  
\rho_3 =1\ .$$

 $\rho_1 $ is instable, $\rho_2 $ is stable  and $\rho_3$ is again instable. 
The domain of attraction of the stable fixpoint $\rho_2$  is   the interval $(\rho_1,1)$.

For $C=F(\rho^{\star})$ we have $\rho_1 = \rho_2 = \rho^{\star}$ and 
all solutions with initial values above $\rho^{\star}$ converge towards $\rho^{\star}$, all other solutions  diverge.

\begin{remark}
In the Lighthill-Whitham case $F(\rho) = \rho(1-\rho)$ we have $$\rho_{1,2} = \frac{1}{2} (1 \mp \sqrt{1-4 C})$$
with $C< \frac{1}{4}$.
For $C=\frac{1}{4}$ we have $\rho_1 = \rho_2 = \frac{1}{2}$.

More explicitly, the layer solution could be determined by solving equation (\ref{layer}). This is   in the Lighthill-Whitham  case a so called Abel
differential equation with constant coefficients. It  can be solved explicitly to obtain  the detailed behavior of the solution in the layer.
	
\end{remark}

\begin{remark}
For the right boundary at $x_R$  a scaling $x \rightarrow \frac{x_R-x}{\epsilon}$ gives 
the layer equations for $x \in[0, \infty)$ as

\begin{align}\label{layerright}
\begin{aligned}
 q &=C\\
 - \partial_x \rho   &= (1-\rho) \frac{F(\rho)-C}{HC}\ .
\end{aligned}
\end{align}
For  $C < F(\rho^{\star})$
the above problem has again three  fix points
$$\rho_{1} \le \rho^{\star}\ ,\ 
\rho_2 = \tau (\rho_1) \ge \rho^{\star}\ ,\ 
\rho_3 =1\ .$$
In this case 
 $\rho_1 $ is stable, $\rho_2 $ is instable and $\rho_3$ is again stable. 
The domain of attraction of the stable fixpoint $\rho_1$  is   $(0,\rho_2)$.
The domain of attraction of the stable fixpoint $\rho_3$  is   $(\rho_2,1)$.

For $C=F(\rho^{\star})$ we have $\rho_1 = \rho_2 = \rho^{\star}$ and 
all solutions with initial values below  $\rho^{\star}$ converge towards $\rho^{\star}$, all other solutions  converge towards $1$.

\end{remark}

\subsection{Macroscopic boundary conditions for the nonlinear kinetic equations}

For the left boundary we prescribe
for the kinetic equation the 2- Riemann invariant $g_2(x_L) = \frac{Hq(x_L)}{(1-\rho(x_L))^H} =  \frac{H f_2(x_L)}{(1-f_1(x_L)-f_2(x_L))^H}$ and for the right boundary the 1-Riemann invariant $g_1(x_R)=f_1(x_R)$.
The boundary conditions for the scalar problem are now derived from the kinetic ones by considering the layer equations with the above boundary conditions at $x=0$ and determining the asymptotic state $\rho_K$, i.e. the solution at infinity of the layer equations.
This state is then used as boundary condition for the scalar equations. The inital value of the scalar equation is in the following denoted by $\rho_B$.

\subsubsection{Left boundary}
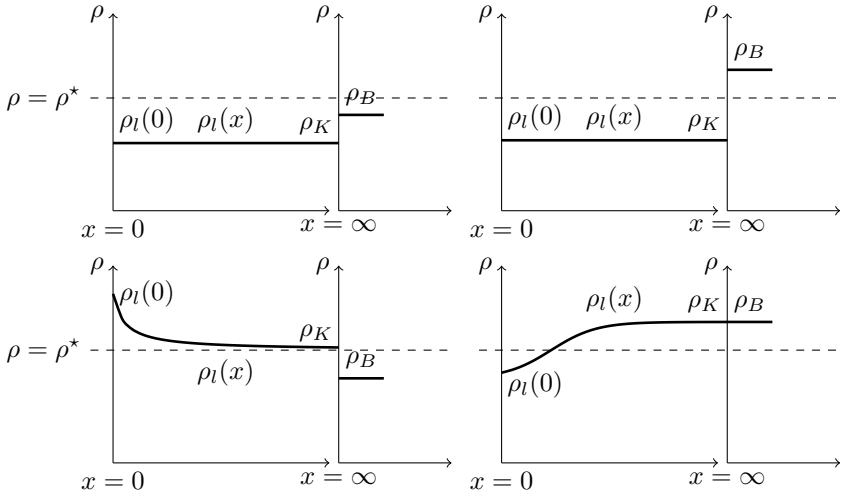
\begin{figure}[h]
	\externaltikz{Layer_cases}{
	\begin{tikzpicture}[scale = .75]
	\def\len{4}
	\def\low{-2}
	
	\def \rhozero {-1}
	\def \rhoB {-0.3}
	\def \rhoK {-0.8}
	\def \rhok {-0.75}
	\def \Ckonst {\rhok*\rhok}
	\def \Ctwo {-1/(2*\rhok)*ln((\rhok+\rhozero)/(\rhok-\rhozero))}

	\node[] (xinf) at (\len,\low){};
	\node[below] at (0,\low){$x=0$};
	\node[below] at (\len,\low){$x=\infty$};
	\draw[dashed] (-0.1*\len,0)--(1.5*\len,0) node[left,pos = 0]{$\rho=\rho^\star$};
	\draw[->] (0,\low)--(xinf);
	\draw[->] (\len,\low)--(1.5*\len,\low);
	\draw[->] (0,\low)--(0,1.5) node[left]{$\rho$};
	\draw[->] (\len,\low)--(\len,1.5) node[left]{$\rho$};
	\draw[line width=1pt] (\len,\rhoB)--(1.2*\len,\rhoB) node[above,pos=0.5] {$\rho_B$};
	\draw[line width=1pt,domain=0.0:\len,smooth,variable=\x,] plot ({\x},{\rhoK}) ;
	\node[above] at (0.15*\len,\rhoK) {$\rho_l(0)$};
	\node[above] at (0.5*\len,\rhoK) {$\rho_l(x)$};
	\node[above] at (\len*0.9,\rhoK) {$\rho_K$};	
	\end{tikzpicture}
	\begin{tikzpicture}[scale = .75]
				\def\len{4}
				\def\low{-2}
				
				\def \rhozero {-1}
				\def \rhoB {0.5}
				\def \rhok {-0.75}
				\def \Ckonst {\rhok*\rhok}
				\def \Ctwo {-1/(2*\rhok)*ln((\rhok+\rhozero)/(\rhok-\rhozero))}

				\node[] (xinf) at (\len,\low){};
				\node[below] at (0,\low){$x=0$};
				\node[below] at (\len,\low){$x=\infty$};
				\draw[dashed] (-0.1*\len,0)--(1.5*\len,0) node[left,pos = 0]{$$};
				\draw[->] (0,\low)--(xinf);
				\draw[->] (\len,\low)--(1.5*\len,\low);
				\draw[->] (0,\low)--(0,1.5) node[left]{$\rho$};
				\draw[->] (\len,\low)--(\len,1.5) node[left]{$\rho$};
				\draw[line width=1pt] (\len,\rhoB)--(1.2*\len,\rhoB) node[above,pos=0.5] {$\rho_B$};
				\draw[line width=1pt,domain=0.0:\len,smooth,variable=\x,] plot ({\x},{\rhok}) ;
				\node[above] at (0.15*\len,\rhok) {$\rho_l(0)$};
				\node[above] at (0.5*\len,\rhok) {$\rho_l(x)$};
				\node[above] at (\len*0.9,\rhok) {$\rho_K$};	
				\end{tikzpicture}
				
	\begin{tikzpicture}[scale = .75]
	\def\len{4}
	\def\low{-2}
	
	\def \rhozero {-1}
	\def \rhoB {-0.5}
	\def \rhok {0}
	\def \Ckonst {\rhok*\rhok}
	\def \Ctwo {-1/(2*\rhok)*ln((\rhok+\rhozero)/(\rhok-\rhozero))}

	\node[] (xinf) at (\len,\low){};
	\node[below] at (0,\low){$x=0$};
	\node[below] at (\len,\low){$x=\infty$};
	\draw[dashed] (-0.1*\len,0)--(1.5*\len,0) node[left,pos = 0]{$\rho = \rho^\star$};
	\draw[->] (0,\low)--(xinf);
	\draw[->] (\len,\low)--(1.5*\len,\low);
	\draw[->] (0,\low)--(0,1.5) node[left]{$\rho$};
	\draw[->] (\len,\low)--(\len,1.5) node[left]{$\rho$};
	\draw[line width=1pt] (\len,\rhoB)--(1.2*\len,\rhoB) node[above,pos=0.5] {$\rho_B$};
	\draw[line width=1pt,domain=0.0:\len,smooth,variable=\x,] plot ({\x},{1/(\x/\len*20-1/(\rhozero))}) ;
	\node[below] at (0.15*\len,1.4) {$\rho_l(0)$};
	\node[below] at (0.5*\len,\rhok) {$\rho_l(x)$};
	\node[above] at (\len*0.9,\rhok) {$\rho_K$};	
	\end{tikzpicture}
		\begin{tikzpicture}[scale = .75]
		\def\len{4}
		\def\low{-2}
		
		\def \rhozero {0.4}
		\def \rhoB {0.5}
		\def \rhok {0.5}
		\def \Ckonst {\rhok*\rhok}
		\def \Ctwo {-1/(2*\rhok)*ln((\rhok+\rhozero)/(\rhok-\rhozero))}

		\node[] (xinf) at (\len,\low){};
		\node[below] at (0,\low){$x=0$};
		\node[below] at (\len,\low){$x=\infty$};
		\draw[dashed] (-0.1*\len,0)--(1.5*\len,0) node[left,pos = 0]{$$};
		\draw[->] (0,\low)--(xinf);
		\draw[->] (\len,\low)--(1.5*\len,\low);
		\draw[->] (0,\low)--(0,1.5) node[left]{$\rho$};
		\draw[->] (\len,\low)--(\len,1.5) node[left]{$\rho$};
		\draw[line width=1pt] (\len,\rhoB)--(1.2*\len,\rhoB) node[above,pos=0.5] {$\rho_B$};
			\draw[line width=1pt,domain=0.0:\len,smooth,variable=\x,] plot ({\x},{-\rhok*tanh(-\rhok*(\x/\len*10+\Ctwo))}) ;
		\node[above] at (0.15*\len,-1.0) {$\rho_l(0)$};
		\node[above] at (0.5*\len,\rhoB) {$\rho_l(x)$};
		\node[above] at (\len*0.9,\rhoB) {$\rho_K$};	
		\end{tikzpicture}
	}	
		
	\caption{Boundary layer and Riemann problem solution for the different cases at the left boundary. First row: case 1 a) and b). Second row case 2 and 3. }
	\label{layerfig}
\end{figure}

Assume for the left boundary  $0 \le g_2(x_L)$ to be known. 
We distinguish three cases. An illustration of the different situations is given in Figure \ref{layerfig}.

\begin{figure}[h]
	\center
	\externaltikz{fundamental_LWR}{
	\begin{tikzpicture}[scale = 6]
	\def \rhobar {0.3}	
	\def \rhostar {0.5}
	\draw[->] (0,0)--(1.2,0) node[below]{$\rho$};
	\draw[->] (0,0)--(0,0.3) node[left]{$F(\rho)$}node at (0.2,0.25) {$F(\rho_1)$};
	\draw[black,line width=1pt,domain=0.0:1,smooth,variable=\x,] plot ({\x},{\x*(1-\x)}) ;
	\draw[dashed] (\rhobar,{\rhobar*(1-\rhobar)})--(\rhobar,0) node[below]{$\rho_1$};
	\draw[dashed] (\rhostar,{\rhostar*(1-\rhostar)})--(\rhostar,0) node[below]{$\rho^*$}node at (0.5,0.3) {$F(\rho^*)$};;
	\draw[dashed] (\rhobar,{\rhobar*(1-\rhobar)})--(1,0)
     node[below]{$1$};	;
	\end{tikzpicture}
	}
	\caption{Fundamental diagram.}
	\label{figfund}
\end{figure}
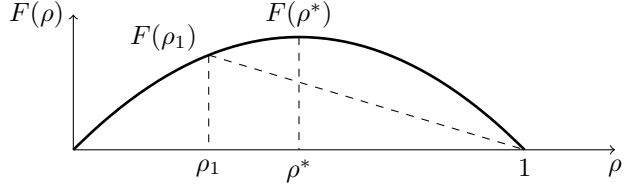

\paragraph{Case 1: ingoing flow}
 $\rho_B < \rho^{\star}$ and $0 \le g_2(x_L) \le \frac{HF(\rho^{\star})}{(1-\rho^{\star})^H}$ or  $\rho_B >\rho^{\star}$ and $g_2(x_L)\le \frac{HF(\tau(\rho_B))}{(1-\tau(\rho_B))^H} $.

The  layer solution is in this case the unstable solution
$$
\rho_l(x) = \rho_1 (C) \le \rho^{\star}\ .
$$
Here,  $0<C< F(\rho^{\star})$ is determined from $g_2(x_L) 
=  \frac{HC}{(1-\rho_1(C))^H } $. First, we determine $\rho_1$ from 
$\frac{HF(\rho_1)}{(1-\rho_1)^H}= g_2(x_L)$. This has a unique solution due to the assumption on $g_2(x_L)$, see Figure \ref{figfund}.
Then,   $C$ is determined from $\rho_1$.
We have  
$\rho_L(0) = \rho_1 \le \rho^{\star}$
and
$$\rhoK=\rho_l(0)\ . $$

In the first case $\rho_l(0) \le \rho^{\star}$ and in the second case $\rho_l(0) \le\tau(\rho_B)$.
In both cases one obtains a wave with positive speed starting at the boundary.
\paragraph{Case 2: transonic flow} $\rho_B <\rho^{\star}$ and $g_2(x_L) \ge  \frac{HF(\rho^{\star})}{(1-\rho^{\star})^H}$.

One chooses $C$ as the maximal possible value $C = F(\rho^{\star})$.
From $g_2(x_L) = \frac{HC}{(1-\rho_l(0))^H}$, which gives $g_2(x_L)= \frac{HF(\rho^{\star})}{(1-\rho_l(0))^H}$, we obtain $\rho_l(0) = 1-\left(\frac{HF(\rho^{\star})}{ g_2(x_L)}\right)^{\frac{1}{H}}$.
This yields $\rho_l(0)\ge \rho^{\star}$. 
The layer solution is no longer constant in space.
Moreover,  $\rhoK=\rho^{\star}$.
In this case one obtains a rarefaction wave with $\rho^\star$ at the boundary.

\paragraph{Case 3: outgoing flow}
$\rho_B >\rho^{\star}$ and 
$g_2(x_L)\ge \frac{HF(\tau(\rho_B))}{(1-\tau(\rho_B))^H}$.

Here,
$$
\rho_l(\infty) = \rho_2(C) = \rho_B
$$
yields $C$ and 
$g_2(x_L)= \frac{HC}{(1-\rho_l(0))^H}= \frac{HF(\rho_2)}{(1-\rho_l(0))^H}= \frac{H F(\rho_B)}{(1-\rho_l(0))^H}$ 
gives
$$
\rho_l(0) = 1- \left(\frac{HF(\rho_B)}{g_2(x_L)}\right)^{\frac{1}{H}}\ge \tau(\rho_B)\ .
$$
Obviously  $\rhoK = \rho_B$. There is no wave starting at the boundary and the layer does not have a constant solution.

\subsubsection{Right boundary}

For the right boundary we prescribe the 1-Riemann invariant $g_1(x_R)=f_1(x_R)= \rho(x_R) -q(x_R)$.

\paragraph{Case 1: ingoing flow}
 $\rho_B > \rho^{\star}$ and $ 1 \ge g_1(x_R) \ge \rho^{\star}- F(\rho^{\star})$ or  $\rho_B <\rho^{\star}$ and $g_1(x_R)\ge \tau(\rho_B)- F(\tau(\rho_B)) $.

The layer solution is 
$$
\rho_r(x) = \rho_2 (C) \ge \rho^{\star}.
$$
Here,  $0<C< F(\rho^{\star})$ is determined from $g_1(x_R) 
=  \rho_2(C) -C $. We determine $\rho_2$ from 
$\rho_2-F(\rho_2)= g_1(x_R)$. This has a unique solution due to the assumption on $g_1(x_R)$.
Then,   $C$ is determined from $\rho_2$.
Moreover, 
$\rho_r(0) = \rho_2 \ge \rho^{\star}$
and
$$\rhoK=\rho_r(0)\ . $$

In the first case $\rho_r(0) \ge \rho^{\star}$ and in the second case $\rho_r(0) \ge\tau(\rho_B)$.
\paragraph{Case 2: transonic flow} $\rho_B >\rho^{\star}$ and $g_1(x_R) \le \rho^{\star}-  F(\rho^{\star})$.

In this case we have $C = F(\rho^{\star})$.
From $g_1(x_R) = \rho_r(0)-C= \rho_r(0)-F(\rho^{\star})$ we obtain $\rho_r(0) = g_1(x_R)+ F(\rho^{\star})$.
This yields $\rho_r(0)\le \rho^{\star}$. 
Moreover,  $\rhoK=\rho^{\star}$.


\paragraph{Case 3: outgoing flow}
$\rho_B <\rho^{\star}$ and 
$g_1(x_R)\le \tau(\rho_B)- F(\tau(\rho_B)) $.

Then,
$$
\rho(\infty) = \rho_1(C) = \rho_B\ .
$$
This yields $C$
and 
$g_1(x_R)= \rho_r(0)-C= \rho_r(0)- F(\rho_1)= \rho_r(0)- F(\rho_B)$ gives
$$
\rho_r(0) = g_1(x_R)+ F(\rho_B)\le \tau(\rho_B)\ .
$$
We have  $\rhoK = \rho_B$.

\subsection{Boundary conditions for the Lighthill Whitham case}
First, the left boundary is considered:
\paragraph{Case 1: ingoing flow}
$\rho_B < \frac{1}{2}$ and $0 \le g_2(x_L) \le H 2^{H-2}$ or  $\rho_B >\frac{1}{2}$ and $g_2(x_L)\le  H (\rho_B)^{1-H} (1-\rho_B)$

The  layer solution is
$$
\rho_l(x) = \rho_1 = \frac{1}{2} (1-\sqrt{1-4 C})\ .
$$
This gives $C= \rho_1 - \rho_1^2 \le \frac{1}{4}$.
Here,  $\rho_1$ is determined from $g_2(x_L) 
=  \frac{H\rho_1}{((1-\rho_1)^{H-1}} $, which has a unique solution due to the assumptions on $g_2(x_L)$.
We note that for $H=1$ we have $\rho_1 = g_2(x_L)$, for $H=2$ we have $\rho_1 = \frac{g_2(x_L)}{2+g_2(x_L)}$ and for  $H=3$ we have
$\rho_1 = \frac{1}{2 g_2(x_L)} (2 g_2(x_L) - \sqrt{12 g_2(x_L)+9}+3)$. This  gives 
$$\rhoK=\rho_l(0)= \rho_1\ . $$


\paragraph{Case 2: transonic flow} $\rho_B <\frac{1}{2}$ and $g_2(x_L) \ge H 2^{H-2}$.

From $C = \frac{1}{4}$ and
 $g_2(x_L) = \frac{HC}{(1-\rho_l(0))^H}= \frac{H}{4(1-\rho_l(0))^H}$ we obtain $\rho_l(0) = 1-\left(\frac{H}{4 g_2(x_L)}\right)^{\frac{1}{H}} \ge  \frac{1}{2}$. 
As before  $\rhoK=\frac{1}{2}$.


\paragraph{Case 3: outgoing flow}
$\rho_B >\frac{1}{2}$ and $g_2(x_L)\ge  H (\rho_B)^{1-H} (1-\rho_B) $.

The layer solution is 
$$
\rho_l(\infty) = \frac{1}{2} (1+\sqrt{1-4 C}) = \rho_B\ .
$$
This gives $C = \rho_B - \rho_B^2$ and 
$g_2(x_L)= \frac{HC}{(1-\rho_l(0))^H}= \frac{ H(\rho_B - \rho_B^2)}{(1-\rho_l(0))^H}$ gives
$$
\rho_l(0) = 1- \left(\frac{H(\rho_B - \rho_B^2)}{g_2(x_L)} \right)^{\frac{1}{H}} \ge 1- \rho_B\ .
$$
Again $\rhoK = \rho_B$.

At the right boundary we have

\paragraph{Case 1: ingoing flow}
$\rho_B > \frac{1}{2}$ and $1\ge g_1(x_R) \ge \frac{1}{4}$ or  $\rho_B <\frac{1}{2}$ and $g_1(x_R)\ge  (1-\rho_B)^2 $.

the layer solution is 
$$
\rho_r(x) = \frac{1}{2} (1+\sqrt{1-4 C})\ .
$$
  $0<C< \frac{1}{4}$ is determined from $g_1(x_R) 
=  \frac{1}{2} (1+\sqrt{1-4 C})  -C  $. This gives
$$
C=\rho_r(0) - (\rho_r(0))^2 = \sqrt{g_1(x_R)}-g_1(x_R) \le \frac{1}{4}
$$
due to the assumption on $g_1(x_R)$.
Moreover,  $g_1(x_R) = (\rho_r(0))^2$ 
and
$$\rho_K=\rho_r(0)\ . $$

\paragraph{Case 2: transonic flow} $\rho_B >\frac{1}{2}$ and $g_1(x_R) \le \frac{1}{4}$.

From $C = \frac{1}{4}$ and 
 $g_1(x_R) = \rho_r(0)- \frac{1}{4}$ we obtain $\rho_r(0) = g_1(x_R)+\frac{1}{4} \le \frac{1}{2}$. 
Moreover,  $\rhoK=\frac{1}{2}$.


\paragraph{Case 3: outgoing flow}
$\rho_B <\frac{1}{2}$ and $g_1(x_R)\le (1-\rho_B)^2$.

We have
$$
\rho_r(\infty) = \frac{1}{2} (1-\sqrt{1-4 C}) = \rho_B\ .
$$
This gives $C = \rho_B - \rho_B^2$ and
with 
$g_1(x_R)= \rho_r(0)-C$ one obtains
$$
\rho_r(0) = g_1(x_R)+ \rho_B - \rho_B^2\le 1- \rho_B\ .
$$
As before $\rho_K = \rho_B$.

\section{Relaxation schemes}
\label{relaxationmethods}

The considerations in the previous sections can be also used to design a relaxation method  based on the nonlinear relaxation system
\ref{cons}. For simplicity we consider the special case $H=1$, where the system is totally linear degenerate.
We refer to \cite{B03,B04,CHN09} for relaxation schemes starting from nonlinear totally linear degenerate relaxation systems.

Before describing the scheme, we note that a relaxation method based on the  linear relaxation system 
 \eqref{dvm0} has different drawbacks. First, choosing  $v_1$, $v_2$ according to  the subcharacteristic condition 
 $v_1 \le F^\prime (\rho) \le v_2 $, which means choosing in general $ v_1 <0 $, yields a convergent scheme. However, for $\epsilon >0$
 positivity of $q$ and the restrictions on $\rho$ and $q$ are not guaranteed, see the discussion on the invariant domains in section \ref{notations}.
 Only for $\epsilon = 0$, i.e. for the relaxed scheme , we obtain a reasonable scheme, which is in this case simply a  Lax-Friedrichs type
 scheme.

On the contrary, choosing $v_1 \ge 0$ the scheme would preserve the positivity of $\rho$. However, the scheme would not work for 
negative wave speeds or $\rho > \rho^\star$ due to the violation of the subcharacteristic condition.
Moreover, the restriction  $\rho \le 1$ is not preserved.
We note again that the  solution proposed in \cite{HPS} is not working properly since the underlying relaxation system is unstable in the sense
of ordinary differential equations for negative wave speeds or $\rho > \rho^\star$.

Our nonlinear relaxation scheme is  given by the following considerations.
Split the system \eqref{cons} into an advection and  a relaxation part.
The advection part is solved with the Godunov method. In the totally linear degenerate case $H=1$ this is easily computed as 
\begin{align*}
	\rho^{n+1/2}_i & = \rho^n_i - \frac{\Delta t}{\Delta x}
	\left(q^n_i\frac{1-\rho^n_{i+1}+q^n_{i+1}}{1-\rho^n_{i}+q^n_{i}}-q^n_{n-1}\frac{1-\rho^n_{i}+q^n_{i}}{1-\rho^n_{i-1}+q^n_{i-1}}\right)
	\\
z^{n+1/2}_i & = z^n_i - \frac{\Delta t}{\Delta x}\left(z^n_i-z^n_{i-1}\right)\	.
\end{align*}
The relaxation part is in the simplest case treated by the implicit Euler method
\begin{align*}
\rho^{n+1}_i & = \rho^{n+1/2}_i \\
z^{n+1}_i & = z^{n+1/2}_i  -\frac{\Delta t}{\epsilon}\left( z^{n+1}_i -\frac{ F(\rho^{n+1}_i) }{(1-\rho^{n+1}_i)}  \right)\ .
\end{align*}
Solving the relaxation ODE for $\epsilon = 0$ we obtain $q^n_i = F(\rho^n_i)$.
Thus the relaxed  scheme reads
\begin{align*}
\rho^{n+1}_i & = \rho^n_i - \frac{\Delta t}{\Delta x}
\left(F(\rho^n_i)\frac{1-\rho^n_{i+1}+F(\rho^n_{i+1})}{1-\rho^n_{i}+F(\rho^n_{i})}-f(\rho^n_{n-1})\frac{1-\rho^n_{i}+F(\rho^n_{i})}{1-\rho^n_{i-1}+f(\rho^n_{i-1})}\right)
\ .
\end{align*}
Note that this is neither the Godunov scheme for the limit-equation nor the Lax-Friedrichs scheme.
The accuracy of the above scheme is intermediate between these two scheme, see the numerical investigation in the next section.
Simple computations show that the scheme is consistent. It is monotone if $f(\rho)+(1-\rho)\ f'(\rho)\geq 0$, which is the 
subcharacteristic condition.
Higher order relaxation methods could be derived as well with the usual procedures.
\begin{remark}
If the kinetic equations with $H>1$ are used for the relaxation method, a nonlinear equation of order $H$ has to be solved to find the intermediate state for the Godunov scheme. In special cases this could be done explicitly. In the general case, one has to solve the algebraic equation numerically.
\end{remark}

\section{A constrained model for $H$ going to $0$ and cluster dynamics}
\label{cluster}

In this section we consider the limit $H\rightarrow 0$ and the case without relaxation term. In this case the influence of the braking term is concentrated at the maximal density. This leads to a cluster dynamic. We refer to  \cite{Ber,Deg} for a similar  investigation for the modified Aw-Rascle model.

Letting $H$ go to $0$ in (\ref{macro0}) one obtains for $\rho<1$ the following simple linear equation

\begin{align}\label{macrohom}
\begin{aligned}
\partial_t \rho + \partial_x q &=0\\
\partial_t q +  \partial_x  q  &=0 \ .
\end{aligned}
\end{align}

This equation would have the invariant domain $ q \le \rho \le q+1$ and $ 0  \le q \le 1$. That means the density could exceed its maximal value $\rho=1$.  However, due to the singularity in the breaking term at $\rho =1 $
the state space of the limit equation is again restricted to $0 \le \rho \le 1$ and $ 0 \le q \le \rho$ as will be discussed in the following.
To find the  dynamics for $H=0$ one has to consider the solution of the Riemann problems for the original system with $H>0$  and let $H$ go  to $0$, compare \cite{Ber,Deg} for such a discussion for the modified Aw-Rascle equation.

First, we remind the reader, that the shock speed $s$   for the 1-wave is given by the following expression, compare  equation (\ref{shockspeed}):
\begin{align}\label{shockspeed00}
\begin{aligned}
 s = 
 \frac{q_L}{(1-\rho_L)^H}\frac{(1-\rho_L)^H - (1-\rho_{R})^H}{\rho_L-\rho_{R}}\ .
\end{aligned}
\end{align}

We consider now a situation with $1-q_L \le \rho_R \le 1$ and $0<q_R < \rho_R+q_L-1$. Outside of this region in state space,  the solution of the Riemann problem is directly described by the solution of the linear formal limit equation (\ref{macrohom}) with  waves with speed $0$ and $1$ and intermediate states given by 
$(\rho_M,q_M) = (\rho_R+ q_L-q_R,q_L)$.

In case the initial values of the Riemann problem are  restricted by $1-q_L \le \rho_R \le 1$ and $0<q_R < \rho_R+q_L-1$, the linear equation would yield a solution with $\rho>1$. In this case  we consider instead  the Riemann problem for the system with $H>0$ and investigate its behaviour as $H$ goes to $0$.

That leads to the following.
The solution of the Riemann problem is given by a 1-shock curve combined with a 2-contact discontinuity.
The intermediate state $\rho_M,q_M$ is given by the intersection of these two curves which gives
\begin{align}\label{intermed1}
\begin{aligned}
 \rho_M+q_R-\rho_R = q_L \frac{(1-\rho_M)^H}{(1-\rho_L)^H}.
\end{aligned}
\end{align}
We do not have to solve that explicitly, just remark that
\begin{align}\label{intermed2}
\begin{aligned}
\frac{1}{q_L} \left((\rho_M+q_R-\rho_R )(1-\rho_L)^H \right)=(1-\rho_M)^H.
\end{aligned}
\end{align}

Using this in the shock speed equation (\ref{shockspeed00}) with $\rho_R = \rho_M$ we obtain

\begin{align}\label{shockspeedb}
\begin{aligned}
 s &= \frac{q_L}{(1-\rho_L)^H}\frac{(1-\rho_L)^H - (1-\rho_{M})^H}{\rho_L-\rho_{M}}\\
 &= \frac{q_L}{(1-\rho_L)^H}\frac{(1-\rho_L)^H -\frac{1}{q_L} \left((\rho_M+q_R-\rho_R )(1-\rho_L)^H \right) }{\rho_L-\rho_{M}}
 \end{aligned}
 \end{align}
 and finally
 \begin{align}\label{shockspeedfinalH}
 \begin{aligned}
 s = \frac{q_L- \rho_M-q_R + \rho_R}{\rho_L-\rho_{M}}\ .
\end{aligned}
\end{align}
Considering (\ref{intermed1}) we remark, that $\rho_M$ converges to $1$ as $H\rightarrow 0$.
The shockspeed is then for $H\rightarrow 0$ (and $\rho_M \rightarrow 1$ ) given by
\begin{align}\label{shockspeedfinal}
 \begin{aligned}
 s = \frac{1-q_L+q_R - \rho_R}{1-\rho_L}\ .
\end{aligned}
\end{align}
Since the  initial values $(\rho_L,q_L)$ and $(\rho_R,q_R)$ under consideration are restricted  by
$1-q_L \le \rho_R \le 1$ and $0 \le q_R \le \rho_R+q_L-1$
we obtain $s \le 0$.

In conclusion,  the solution of the constrained model for $H=0$  is given by the solution of the linear model as long as the resulting intermediate states are in $0\le \rho \le 1$, $0 \le q \le \rho$ that means for $\rho_R-q_R < 1-q_L$. In this case we have waves with speed $0$ and $1$ and intermediate states given by $(\rho_M,q_M) = (\rho_R+ q_L-q_R,q_L)$.  
In the other cases with $\rho_R-q_R > 1-q_L$, we have a solution with an intermediate state given by $(\rho_M,q_M) = (1, 1+q_R-\rho_R)$. The solution   is a combination of a shock solution with speed $s= \frac{1-q_L+q_R - \rho_R}{1-\rho_L}< 0 $ and a contact discontinuity with speed $s=1$.

We refer again to  \cite{Deg,Ber} for similar investigations for   the modified Aw-Rascle model and for further references on constrained models.

\begin{remark}
The cases are distinguished by determining whether $\rho_R -q_R < 1-q_L$ or $\rho_R -q_R > 1-q_L$.
With $f_1 =\rho-q$, the number of stopped vehicles, this condition could be interpreted as follows: the number of stopped cars on the right is smaller (respectively larger) than the  number of stopped cars on the left  for a left state with maximal density  and a flux equal to $q_L$.
\end{remark}

\begin{remark}
If $\rho_R=1$, then the intermediate state is $(\rho_M,q_M) = (1, q_R)$ and the resulting shock speed is $s= \frac{q_R- q_L}{1-\rho_L}< 0 $. 
\end{remark}

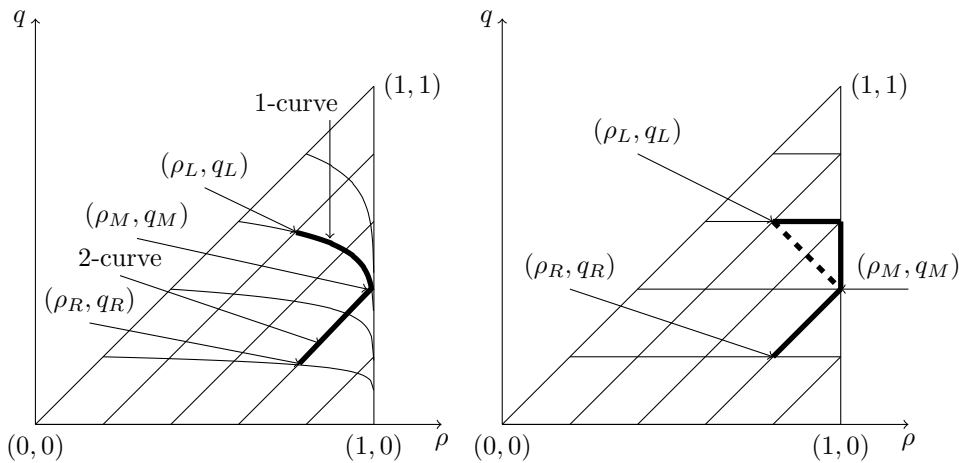
\begin{figure}[h]
	\externaltikz{Statespace_H0}{
	\begin{tikzpicture}[scale = 4.5]
 	\def\dr{0.2}
 	\def\ds{0.1}
 	\node[below] at (0,0) {$(0,0)$};
 	\node[below] at (1,0) {$(1,0)$};
 	\node[right] at (1,1) {$(1,1)$};
 	\draw[->](0,0)--(1.2,0) node[below]{$\rho$};
 	\draw[->](0,0)--(0,1.2) node[left]{$q$};
 	\draw[->](0.25,0.44)--(0.84,0.24) node[above]  at (0.25,0.44){$2$-curve};
 	\draw[->](0.87,0.9)--(0.87,0.55) node[above] at (0.77,0.9){$1$-curve};
 	\draw[->](0.5,0.7)--(0.77,0.57) node[above] at (0.5,0.7){$(\rho_L,q_L)$};
 	 	\draw[->](0.16,0.3)--(0.78,0.18) node[above] at (0.16,0.3){$(\rho_R,q_R)$};
 	 	\draw[->](0.3,0.55)--(0.98,0.4) node[above] at (0.3,0.55){$(\rho_M,q_M)$};
 	 	\draw[line width=2pt] (0.99,0.4)--(0.78,0.18);
 	\draw(1,0)--(1,1);
 	\foreach \r in {0, \dr,...,1}
 	\draw(\r,0)--(1,1-\r);
 	\foreach \r in {0, \dr,...,0.9}
 \draw[domain=\r:1,smooth,variable=\x,samples=100] plot ({\x},{\r*(1-\x)^\ds/(1-\r)^\ds});
 \draw[line width=2pt,domain=0.77:0.994,smooth,variable=\x,samples=100] plot ({\x},{0.6*(1-\x)^\ds/(1-0.6)^\ds});
 	\end{tikzpicture}
 	\begin{tikzpicture}[scale = 4.5]
 	\def\dr{0.2}
 	\node[below] at (0,0) {$(0,0)$};
 	\node[below] at (1,0) {$(1,0)$};
 	\node[right] at (1,1) {$(1,1)$};
 	\draw[->](0,0)--(1.2,0) node[below]{$\rho$};
 	\draw[->](0,0)--(0,1.2) node[left]{$q$};
 	\draw[->](0.4,0.8)--(0.8,0.6) node[above] at (0.4,0.8){$(\rho_L,q_L)$};
 	\draw[->](0.2,0.4)--(0.8,0.2) node[above] at (0.2,0.4){$(\rho_R,q_R)$};
 	\draw[->](1.2,0.4)--(1.0,0.4) node[above] at (1.2,0.4){$(\rho_M,q_M)$};
 	\draw[line width=2pt] (0.8,0.6)--(1.0,0.6);
 	\draw[line width=2pt,dashed] (0.8,0.6)--(1.0,0.4);
 	\draw[line width=2pt] (1.0,0.6)--(1.0,0.4);
 	\draw[line width=2pt] (1.0,0.4)--(0.8,0.2);
 	\draw(1,0)--(1,1);
 	\foreach \r in {0, \dr,...,1}
 	\draw(\r,0)--(1,1-\r);
 	\foreach \r in {0, \dr,...,1}
 	\draw(\r,\r)--(1,\r);
 	\end{tikzpicture}
	 }
 
 	\caption{State space and Riemann problem for  $(\rho,q)$  variables for $H$ small ($H=0.1$) and $H=0$.
 	Situation with $\rho_R-q_R > 1-q_L.$}
 	\label{state}
 \end{figure}

\section{Numerical results}
\label{numericalresults}
In this section we discuss four topics. First, the 
 solution of Riemann problems for the kinetic problem with different $H$ and different $\epsilon$  are discussed.
 Second we investigate the boundary value problem and compare kinetic and limit equations with boundary conditions from Section \ref{boundaryconditions}.
 Third, we present an investigation of the relaxation scheme   for the Lightill-Whitham problem.
 Last, the constrained equations for $H=0$ are investigated
 and compared to the solutions for small $H$.

If not otherwise stated we use the relaxation scheme for the kinetic equations, i.e. the  Godunov scheme for the advection part in conservative form
and the implicit Euler method for the right hand side. For $H=1$ the formulas stated in the last section are used. For the other cases
the resulting algebraic equations of order $H$ are solved using Bisection method. 

The solutions of the limit equation are given as exact solutions.
The number of cells is $1000$ and  if not specified differently we use for the kinetic problem $H=1$ and $\epsilon=0.1$.
The CFL condition is chosen with a CFL number $1$.

\subsection{Numerical solution of Riemann problems for the kinetic equation with different $H$   and different $\epsilon$}

We consider two different Riemann Problems,
first $\rho_L=0.99$ and  $\rho_R=0$ with $q\equiv 0$ and 
second $\rho_L=0.3$ and  $\rho_R=0.99$ with $q\equiv 0$.
In the first example the LWR-solution is a rarefaction wave, in the second case it is a shock wave.
Note that, if for the second example the left and right states in the kinetic equation are in equilibrium, then, the speed of the left going shock wave coincides for any $\epsilon$ with the shock speed in the LWR model:
from the first equation of the Rankine-Hugoniot condition \eqref{locus} we obtain for the left going wave
\begin{align*}
	s = \frac{q_R-q_L}{\rho_R-\rho_L}\ .
\end{align*}
If $q$ is as in the LWR model also $s$ is identical.
		\begin{figure}[h!]
			\externaltikz{Rare_eps_1}{
				\begin{tikzpicture}[scale=0.6]
				\begin{axis}[ylabel = $\rho$,xlabel =  $x$,
				legend style = {at={(0.5,1)},xshift=0.0cm,yshift=0.1cm,anchor=south},
				legend columns= 4,
				]
				\addplot[color = red,thick] file{Data/rho_Rare_eps05.txt};
				\addlegendentry{$\epsilon = 0.5$}
				\addplot[color = blue,thick] file{Data/rho_Rare_eps01.txt};
				\addlegendentry{$\epsilon = 0.1$}
				\addplot[color = green,thick] file{Data/rho_Rare_eps001.txt};
				\addlegendentry{$\epsilon = 0.01$}
				\addplot[color = orange,thick] file{Data/rho_Rare_eps0001.txt};
				\addlegendentry{$\epsilon = 0.001$}
				\end{axis}
				\end{tikzpicture}
			}
			\externaltikz{shock_eps_1}{
				\begin{tikzpicture}[scale=0.6]
				\begin{axis}[ylabel = $\rho$,xlabel =  $x$,
				legend style = {at={(0.5,1)},xshift=0.2cm,yshift=0.1cm,anchor=south},
				legend columns= 4,
				]
				\addplot[color = red,thick] file{Data/rho_shock_eps05.txt};
				\addplot[color = blue,thick] file{Data/rho_shock_eps01.txt};
				\addplot[color = green,thick] file{Data/rho_shock_eps001.txt};
				\addplot[color = orange,thick] file{Data/rho_shock_eps0001.txt};
				\end{axis}
				\end{tikzpicture}
			}
			\caption{Numerical solutions for different $\epsilon$ at $t=0.4$.}
			\label{epsilon}
		\end{figure}
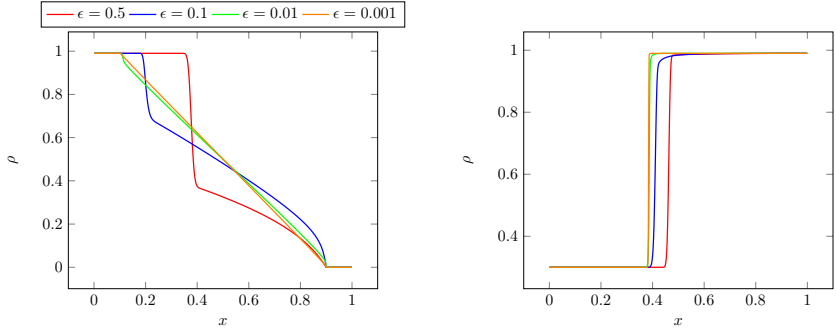
Figure \ref{epsilon} shows the two examples for different values of $\epsilon$.
		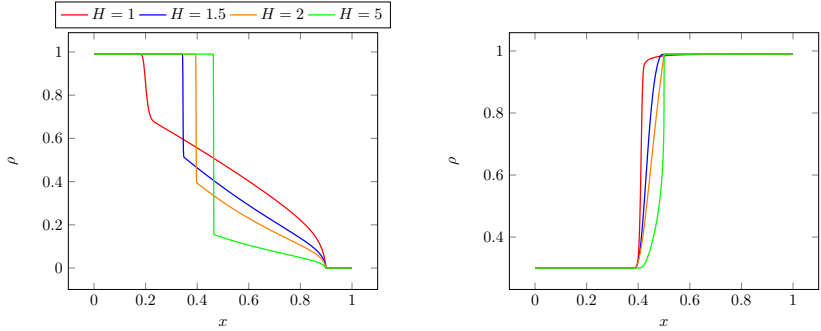
\begin{figure}[h!]
			\externaltikz{Rare_H_1}{
				\begin{tikzpicture}[scale=0.6]
				\begin{axis}[ylabel = $\rho$,xlabel =  $x$,
				legend style = {at={(0.5,1)},xshift=0.0cm,yshift=0.1cm,anchor=south},
				legend columns= 4,
				]
				\addplot[color = red,thick] file{Data/rho_Rare_eps01.txt};
				\addlegendentry{$H = 1$}
				\addplot[color = blue,thick] file{Data/rho_Rare_H15.txt};
				\addlegendentry{$H = 1.5$}
				\addplot[color = orange,thick] file{Data/rho_Rare_H2.txt};
				\addlegendentry{$H = 2$}
				\addplot[color = green,thick] file{Data/rho_Rare_H5.txt};
				\addlegendentry{$H = 5$}
				\end{axis}
				\end{tikzpicture}
			}
			\externaltikz{shock_H_1}{
				\begin{tikzpicture}[scale=0.6]
				\begin{axis}[ylabel = $\rho$,xlabel =  $x$,
				legend style = {at={(0.5,1)},xshift=0.2cm,yshift=0.1cm,anchor=south},
				legend columns= 4,
				]
				\addplot[color = red,thick] file{Data/rho_shock_eps01.txt};
				\addplot[color = blue,thick] file{Data/rho_shock_H15.txt};
				\addplot[color = orange,thick] file{Data/rho_shock_H2.txt};
				\addplot[color = green,thick] file{Data/rho_shock_H5.txt};
				\end{axis}
				\end{tikzpicture}
			}
			\caption{Numerical solutions for different $H$ at $t=0.4$.}
			\label{H}
		\end{figure}
Figure \ref{H} shows the two examples for different values of $H$.

\subsection{Comparison of BVP for kinetic and macroscopic equation}

\begin{figure}[h!]
	\externaltikz{Rare_Layer_1}{
		\begin{tikzpicture}[scale=0.6]
		\begin{axis}[ylabel = $\rho$,xlabel =  $x$,
		legend style = {at={(0.5,1)},xshift=0.0cm,yshift=0.1cm,anchor=south},
		legend columns= 4,
		]
		\addplot[black,thick] coordinates {(0,0.9) (0.18,0.9) (0.74,0.2) (1,0.2)};
		\addlegendentry{LWR}
		\addplot[color = red,thick] file{Data/rho_Layer_eps01.txt};
		\addlegendentry{$\epsilon = 0.1$}
		\addplot[color = blue,thick] file{Data/rho_Layer_eps001.txt};
		\addlegendentry{$\epsilon = 0.01$}
		\addplot[color = green,thick] file{Data/rho_Layer_eps0001.txt};
		\addlegendentry{$\epsilon = 0.001$}
		\end{axis}
		\end{tikzpicture}
	}
	\externaltikz{Rare_Layer_2}{
		\begin{tikzpicture}[scale=0.6]
		\begin{axis}[ylabel = $\rho$,xlabel =  $x$,
		legend style = {at={(0.5,1)},xshift=0.0cm,yshift=0.1cm,anchor=south},
		legend columns= 4,
		]
		\addplot[black,thick] coordinates {(0,0.5) (0.24,0.2) (0.46,0.2)  (0.46,0.9) (1,0.9)};
		\addplot[color = red,thick] file{Data/rho_Layer_test2_eps01.txt};
		\addplot[color = blue,thick] file{Data/rho_Layer_test2_eps001.txt};
		\addplot[color = green,thick] file{Data/rho_Layer_test2_eps0001.txt};
		\end{axis}
		\end{tikzpicture}
	}
	\caption{Layer solutions for different $\epsilon$.}
	\label{figbound}
\end{figure}
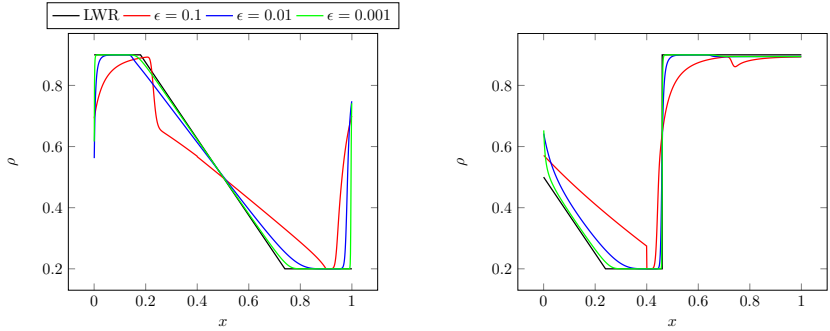

We consider $H=1$. Figure \ref{figbound} shows the solution of a boundary value problem with layers at both boundaries.
In the left picture we have a situation with outgoing flow at the left boundary,
$\rho_B > \frac{1}{2}$  and  $g_2(x_L = 0)$ is chosen such that $\rho(0)\ge 1- \rho_B$.  At the right boundary we have again outgoing flow with $\rho_B < \frac{1}{2}$ and $g_1(x_R=1)$ is chosen such that $\rho(1)\le 1- \rho_B$.

In the right picture the inner states are $\rho(x) = 0.2$ for $x<0.5$ and $\rho(x) = 0.9$ for $x\geq 0.5$.
At the left boundary there is a transsonic flow with $g_2(0) = 0.75$ and with $g_1(1) = 0.8$ we have an ingoing flow at the right boundary.
The figure shows the transonic layer developed at the left boundary.

\subsection{Comparison of numerical schemes} 

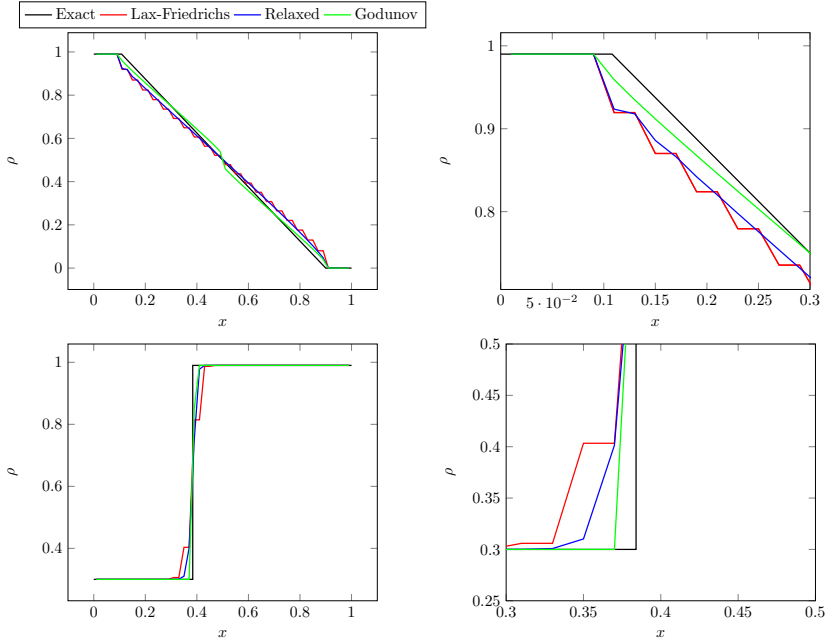
\begin{figure}[h!]
			\externaltikz{Rare_relaxed_1}{
				\begin{tikzpicture}[scale=0.6]
				\begin{axis}[ylabel = $\rho$,xlabel =  $x$,
				legend style = {at={(0.5,1)},xshift=0.0cm,yshift=0.1cm,anchor=south},
				legend columns= 4,
				]
				\addplot[black,thick] coordinates {(0,0.99) (0.108,0.99) (0.9,0.0) (1,0.0)};
				\addlegendentry{Exact}
				\addplot[color = red,thick] file{Data/rho_Rare_eps0_LF.txt};
				\addlegendentry{Lax-Friedrichs}
				\addplot[color = blue,thick] file{Data/rho_Rare_eps0_Relaxed.txt};
				\addlegendentry{Relaxed}
				\addplot[color = green,thick] file{Data/rho_Rare_eps0_Godunov.txt};
				\addlegendentry{Godunov}
				\end{axis}
				\end{tikzpicture}
				\begin{tikzpicture}[scale=0.6]
					\begin{axis}[
					xmin=0.0,
						xmax=0.3,
						xlabel=$x$,
						ylabel=$\rho$,
						]
				\addplot[color = red,thick] file{Data/rho_Rare_eps0_LF.txt};
				\addplot[black,thick] coordinates {(0,0.99) (0.108,0.99) (0.9,0.0) (1,0.0)};
				\addplot[color = red,thick] file{Data/rho_Rare_eps0_LF.txt};
				\addplot[color = blue,thick] file{Data/rho_Rare_eps0_Relaxed.txt};
				\addplot[color = green,thick] file{Data/rho_Rare_eps0_Godunov.txt};
					\end{axis}
				\end{tikzpicture}
			}

			\externaltikz{shock_relaxed_1}{
				\begin{tikzpicture}[scale=0.6]
				\begin{axis}[ylabel = $\rho$,xlabel =  $x$,
				legend style = {at={(0.5,1)},xshift=0.2cm,yshift=0.1cm,anchor=south},
				legend columns= 4,
				]
				\addplot[black,thick] coordinates {(0,0.3) (0.384,0.3) (0.384,0.99) (1,0.99)};
				\addplot[color = red,thick] file{Data/rho_shock_eps0_LF.txt};
				\addplot[color = blue,thick] file{Data/rho_shock_eps0_Relaxed.txt};
				\addplot[color = green,thick] file{Data/rho_shock_eps0_Godunov.txt};
				\end{axis}
				\end{tikzpicture}
				\hspace{0.6cm}  
				\begin{tikzpicture}[scale=0.6]
				\begin{axis}[ylabel = $\rho$,xlabel =  $x$,xmin = 0.3, xmax = 0.5,
				,ymin = 0.25, ymax = 0.5,
				]
				\addplot[black,thick] coordinates {(0,0.3) (0.384,0.3) (0.384,0.99) (1,0.99)};
				\addplot[color = red,thick] file{Data/rho_shock_eps0_LF.txt};
				\addplot[color = blue,thick] file{Data/rho_shock_eps0_Relaxed.txt};
				\addplot[color = green,thick] file{Data/rho_shock_eps0_Godunov.txt};
				\end{axis}
				\end{tikzpicture}
			}
			\label{figrelax}
			\caption{Comparison of Lax-Friedrichs, Godunov and relaxed scheme. The left column shows the solution on the full domain.
			The right column shows zooms at the points of non-smoothness of the solution.}
	\end{figure}

Figure \ref{figrelax} shows a comparison of the Lax-Friedrichs, the Godunov and the relaxed scheme (using $H=1$) from Section \ref{relaxationmethods}.
Two Riemann problems, a rarefaction wave and a shock wave for the Lighthill-Whitham equations are investigated. 
The figure shows that, comparing the two central schemes, the relaxed scheme is more accurate than the Lax-Friedrichs scheme.

\subsection{Cluster dynamic for the constrained equations with $H=0$}

In this section we investigate the limit as $H \rightarrow 0$ numerically and compare the solutions to the 
constrained limit equation for $H=0$ given by the solution of the Riemann problem discussed in Section \ref{cluster}. We consider the two cases $\rho_R-q_R> 1-q_L$ and $\rho_R-q_R< 1-q_L$ with solutions given by the solution of the linear problem (\ref{macrohom}) and solutions given by the discussion in Figure
\ref{state}.We consider for case 1
\begin{align*}
	\rho_L = 0.7\ ,\  q_L=0.7\ ,\
	\rho_R = 0.7\ ,\ & q_R=0.2
\end{align*}
and for case 2 the same values, except $q_L=0.3$.
In both cases a convergence towards the limit solution can be observed.

\begin{figure}
			\externaltikz{Cluster_H}{
				\begin{tikzpicture}[scale=0.6]
				\begin{axis}[ylabel = $\rho$,xlabel =  $x$,
				legend style = {at={(0.5,1)},xshift=0.2cm,yshift=0.1cm,anchor=south},
				legend columns= 4,
				]
				\addplot[black,thick] coordinates {(0,0.7) ({0.5-2/3*0.2},0.7) ({0.5-2/3*0.2},1) (0.7,1) (0.7,0.7) (1,0.7)};
				\addlegendentry{$H = 0$}
				\addplot[color = red,thick] file{Data/rho_Cluster_H1.txt};
				\addlegendentry{$H = 1$}
				\addplot[color = blue,thick] file{Data/rho_Cluster_H05.txt};
				\addlegendentry{$H = 0.5$}
				\addplot[color = green,thick] file{Data/rho_Cluster_H01.txt};
				\addlegendentry{$H = 0.1$}
				\end{axis}
				\end{tikzpicture}
				\begin{tikzpicture}[scale=0.6]
				\begin{axis}[ylabel = $\rho$,xlabel =  $x$,
				legend style = {at={(0.5,1)},xshift=0.2cm,yshift=0.1cm,anchor=south},
				legend columns= 4,
				]
				\addplot[black,thick] coordinates {(0,0.7) ({0.5},0.7) ({0.5},0.8) (0.7,0.8) (0.7,0.7) (1,0.7)};
				\addplot[color = red,thick] file{Data/rho_Cluster_H1_test6.txt};
				\addplot[color = blue,thick] file{Data/rho_Cluster_H05_test6.txt};
				\addplot[color = green,thick] file{Data/rho_Cluster_H01_test6.txt};
				\end{axis}
				\end{tikzpicture}
			}
			\label{figH}
\caption{Solution for different values of $H>0$ and for $H=0$. On the left:
$\rho_R-q_R> 1-q_L$. On the right: $\rho_R-q_R< 1-q_L$. }
\end{figure}
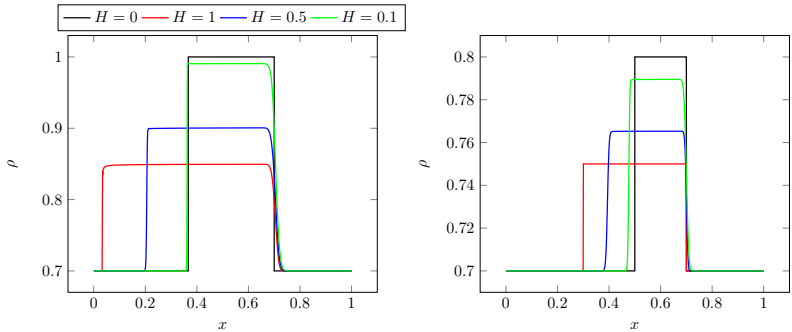

\section{Conclusions}

The paper presents a new nonlinear discrete velocity model for traffic flow having the correct relaxation limit and  having the correct invariant domain for traffic flow modeling. Compared to classical kinetic discrete velocity models it avoids the problems connected with 
the positivity of the velocities and the subcharacteristic condition. In contrast, the hyperbolic part is nonlinear, but relatively simple,
being a totally linear degenerate hyperbolic problem with a simple structure of the integral curves. We have discussed relations to the Aw-Rascle model. Moreover, we have discussed  boundary conditions for the limit equations derived from the relaxation model,we have 
investigated the cluster dynamics of the model for vanishing braking distance and we have suggested a relaxation scheme build on the kinetic discrete velocity model. Numerical results illustrate the behaviour of the solutions for various situations.

\section*{ Acknowledgment}
This research was supported by the German Research Foundation DFG through the SPP 1962.

\end{document}